# Exploring Collaborative and Multidisciplinary Aircraft Optimization through the AGILE Academy Challenge – A case study for an aircraft auxiliary solar power system

*A. Jeyaraj, F. Sanchez, P. Earnest, E. Murugesan, R. Priem, S. Liscouet-Hanke*

**Abstract**

Reduction in aircraft emission is a main driver for the development of more efficient aircraft and enabling technologies are reaching operational maturity. Aircraft manufacturers need an efficient product development process to capture these emergent technologies and develop new aircraft concepts in order to stay competitive. Presently, the aircraft development process is cross organizational, and harnesses distributed, heterogeneous knowledge and expertise. Large scale multidisciplinary studies, involving disciplinary experts and specific tools are required to evaluate different aircraft concepts. These Multi-Disciplinary Aircraft Optimization (MDAO) processes are difficult to deploy as they involve cross-organizational collaboration and harmonization of processes, tools and even vocabulary. Moreover, difficulties in collaborative decision making, reconfiguration and integration of new requirements and competencies often precludes the development of an optimal solution within the available time. To address these challenges, the AGILE project is an effort within the European Union funded Horizon 2020 project to reduce aircraft development time by developing tools and processes that enable efficient, collaborative aircraft design.

This paper presents the work performed as part of the AGILE academy challenge where students were tasked with developing and solving an aircraft MDAO problem using the AGILE toolchain. Each team consisted of students from various universities around the globe and had expertise in multiple design domains. An MDAO study is presented that utilizes the AGILE toolchain to investigate the feasibility of implementing an auxiliary solar power system on a baseline aircraft. The steps performed are: (1) definition of a multi-disciplinary design problem, (2) development of collaborative workflow and (3) optimization using surrogate models. Through the case study, a novel technology concept is investigated and the efficacy of the AGILE toolchain in facilitating a MDAO is analyzed.

# 1 Introduction

Today, the development of safe, environmentally friendly and more efficient aircraft with shorter time to market is required to be successful in a highly competitive market. This can be achieved by integrating new technologies, such as more electrical systems, into future aircraft in an efficient manner. However, the multidisciplinary nature of these technologies presents a challenge to their successful adoption. The current aircraft development process does not incorporate efficient approaches to perform multidisciplinary studies of different aircraft configurations and technologies. Therefore, the aircraft development process requires tools and methodologies to deploy multidisciplinary aircraft optimization strategies to develop the next generation of aircraft. This was the main the main objective of the AGILE project [1]. To this effect, the first objective of AGILE was the structured development of advanced multidisciplinary optimization techniques and their integration to reduce the convergence time in aircraft optimization. Today's advanced analysis codes and software tools are mostly discipline-specific and well understood by disciplinary experts. However; the operation of the system of tools and the interpretation of the results are additional challenges in the collaboration between the disciplinary specialists, and the aircraft generalists. Therefore, the second objective of AGILE is the structured development of processes and techniques for efficient multisite collaboration in the overall design teams [2]. Mastering complex systems depends on the exploitation of knowledge. Besides the interaction of experts, the efficient handling of data, information and knowledge using information technologies shows enormous potential. Thus; the third objective of AGILE is the structured development of knowledge enabled information technologies to support interdisciplinary design campaigns [3]. In order to test the real world validation of the technologies developed within AGILE, the AGILE academy initiative is proposed

The AGILE academy initiative consisted in a series of workshops carried out in collaboration with academic institutions around the world [4]. The objective was to make available all the technologies developed by the AGILE consortium to multiple teams of students who had to work collaboratively on a design task. The AGILE academy challenge represented the second phase of this initiative where multiple teams of students used the AGILE tool chain to solve an aircraft design challenge.
First, this paper describes the AGILE academy challenge and the context of the competition. The aircraft use case is also introduced. Then, the AGILE tools provided by the consortium are introduced. The third part of this paper presents with the case study completed by the team. Finally, discussions and lessons learned are given.

# 2 AGILE academy challenge
## 2.1 Context
The AGILE Challenge, as part of the second cycle of the AGILE Academy initiative, established multiple distributed design teams of students who collaboratively worked and competed with other teams to perform multiple design tasks using the AGILE technologies. Three aircraft design teams were established with students from registered organizations around the world. Each design team was formed by at least 3 different organizations, assembled to guarantee complementary skills in each team. The author's team members are given in Table 1.

Table 1: Team 3 members

| Institution | Country | Member(s) |
|---|---|---|
| University of Pisa | Italy | 2 |
| ISAE SUPAERO | France | 1 |
| ONERA | France | 1 |
| University of Michigan | United States | 1 |
| Concordia University | Canada | 4 |

Each design team was assigned and supported by a main point of contact from the AGILE Partners. The AGILE tool chain was presented to the student teams as an AGILE open source MDO test suite. The distributed collaborative MDO teams were challenged to design an aircraft based on a given baseline requirements set, using the skills and design competencies available within each team. . The aircraft use case is as illustrated in Figure 1.

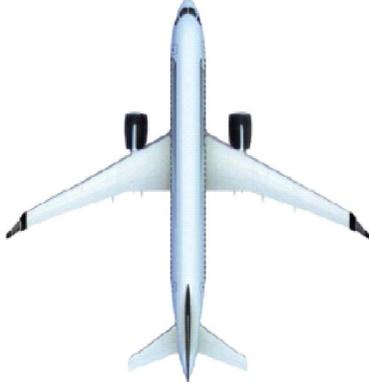

| Conventional Large Regional Jet | |
|---|---|
| Range (102 kg/pax) | 5556 km |
| Design payload | 16329 kg |
| PAX | 10 pax @ 102 kg |
| MLW (% MTOW) | 90 % |
| Cruise Mach (LRC) | 0.78 |
| Initial Cruise Altitude (ICA) | 36000 ft |
| TOFL (ISA, SL, MTOW) | 1900 m |
| LFL (ISA, SL, MTOW) | 1500 m |
| Engine | Turbofan high bypass |

Figure 1: AGILE aircraft use case, adapted from [4]

## 2.2 AGILE tools

The AGILE consortium developed a set of tools to enhance and support collaborative MDA/MDO studies. Figure 2 shows an overview of all the tools and their scope of application within an MDO process; ranging from MDO project organizational tool to remote workflow execution tool. The following sub-sections briefly describe the application of each tool and how they relate to each other.

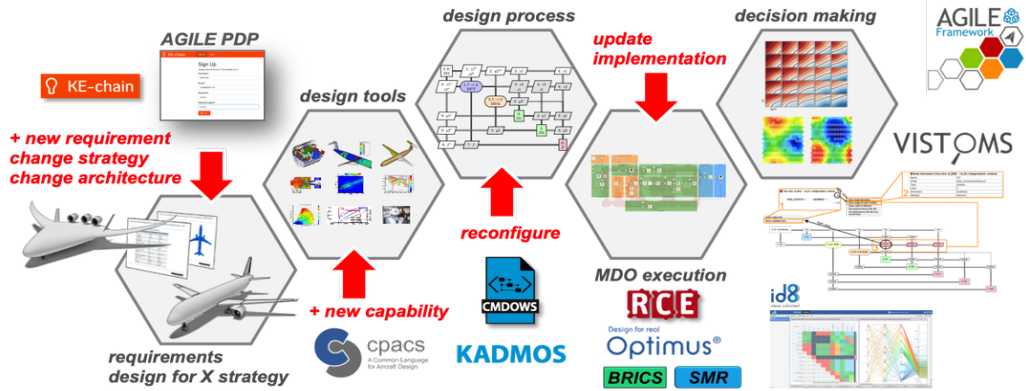

*Figure 2: AGILE tools overview, from [1]*

### 2.2.1 KE-Chain

The AGILE tool package made available for this MDAO project was KE-Chain; an online platform that provided the user with inbuilt functionality to perform MDO activities such as managing project operations, tracking the activities of the different tool operators, structuring the primary design workflow, generating visual representations of the workflow, and producing a downloadable workflow file that can be implemented in a remote workflow executer.

KE-Chain can be perceived as a one-stop-shop, that offers the users all the functionalities vital for an active and collaborate MDO platform, by allowing tool owners to define their tool functionalities and limitations, define user roles, and record each update/activity for traceability. This helps in monitoring the advancement/progress of the project.

### 2.2.2 CPACS: Common Parametric Aircraft Configuration Schema

CPACS is a common language used for the representation and the modeling of parametric aeronautical systems and components, through an intuitive data definition schema. It enables effective exchange of information between the tools and it is therefore a driver for the execution of large-scale collaborative MDO processes. This schema can be structures based aircraft functionalities or components, but all collaborative entities need to follow the exact same data definition strategy in order to exchange data files.

### 2.2.3 CMDOWS: Common MDO Workflow Schema

CMDOWS is an open-source, XML-based storage format for MDO systems. The schema describes a data structure that can be used to store any MDO system at different stages throughout the set-up phase: tool repository, MDO problem, MDO process. In addition, MDO workflow blueprints stored as CMDOWS files can be automatically instantiated as executable workflows in different platforms: RCE, Optimus.

### 2.2.4 KADMOS: Knowledge and graph-based Agile Design with Multidisciplinary Optimization System

KADMOS is an open-source Python-based package that can be used to formulate, inspect and manipulate MDO systems. The package is based on a graph-theoretic formalization of MDO systems. This formalization enables the description of these systems throughout the different stages of the set-up phase in any collaborative project: tool repository, MDO problem, MDO process.

### 2.2.5 VISTOMS: VISualization TOol for MDO Systems

VISTOMS is the graphical user interface that allows the easy setup, inspection and modification of MDO systems. KADMOS and VISTOMS are combined in a free web service called "MDO System Interface". Therefore, all the obtained MDO systems data can be stored and downloaded in the CMDOWS format.

### 2.2.6 RCE: Remote Component Environment

RCE is an open source distributed and workflow-driven integration environment. It is used by design teams to "compose and execute" simulation-based MDO processes of complex systems (e.g., aircraft) by integrating their own design and simulation tools. It can automatically generate MDO workflows stored as CMDOWS format and execute cross-organizational MDO processes.

## 3 Case study: Aircraft auxiliary solar power system

### 3.1 Design problem description

The design problem was to determine the impact of implementing an auxiliary solar power system on the baseline aircraft while keeping certain performance parameters like Take-Off Field length and range the same. The idea is to use the power generated by the SPS to supplement some of the system electrical loads thus resulting in reduced engine power offtake. Additionally, Auxiliary Power Unit (APU) usage on ground is also supplemented by solar power resulting in an overall reduction in fuel consumption which helps reduce the total fuel weight (Figure 3). Therefore, the optimization activities target a maximization of the fuel weight saved while subject to aircraft geometrical constraints to keep it within the same class as the AGILE baseline aircraft.

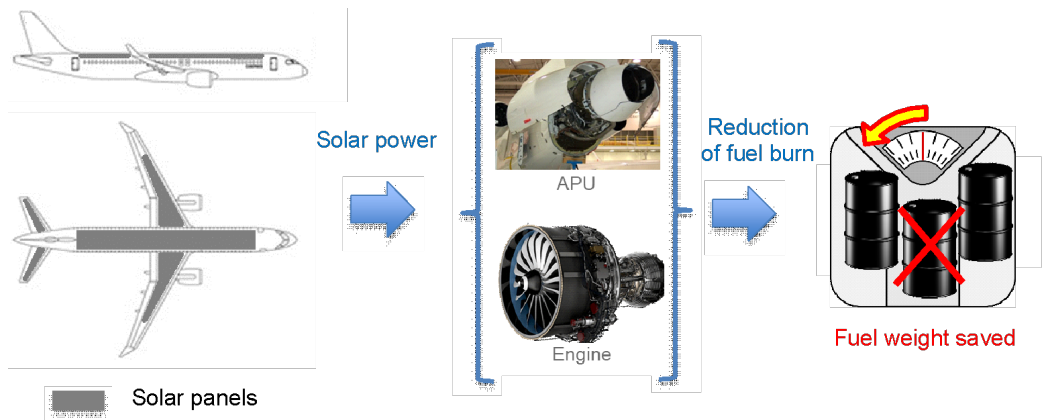

*Figure 3: Design problem of the aircraft auxiliary solar power system*

### 3.2 Task A: MDO formalization

Task A deals with the definition and formalization of an MDO workflow to design the considered aircraft. The preliminary step in building the workflow was to determine which tools will be used and what their interactions will be, making a catalogue of available competences (Figure 4). For the considered design problem, initial studies revealed several tool gaps related to the propulsion competency dealing mainly with the assessment of fuel savings due to reduced engine power offtake. Furthermore, the open-source aircraft sizing tools such as VAMPzero and SUAVE recommended by the AGILE consortium were not suitable for the design of an aircraft auxiliary solar power system. Therefore, ISAE SUPAERO and Concordia University developed two tools to address these requirements (Figure 4).

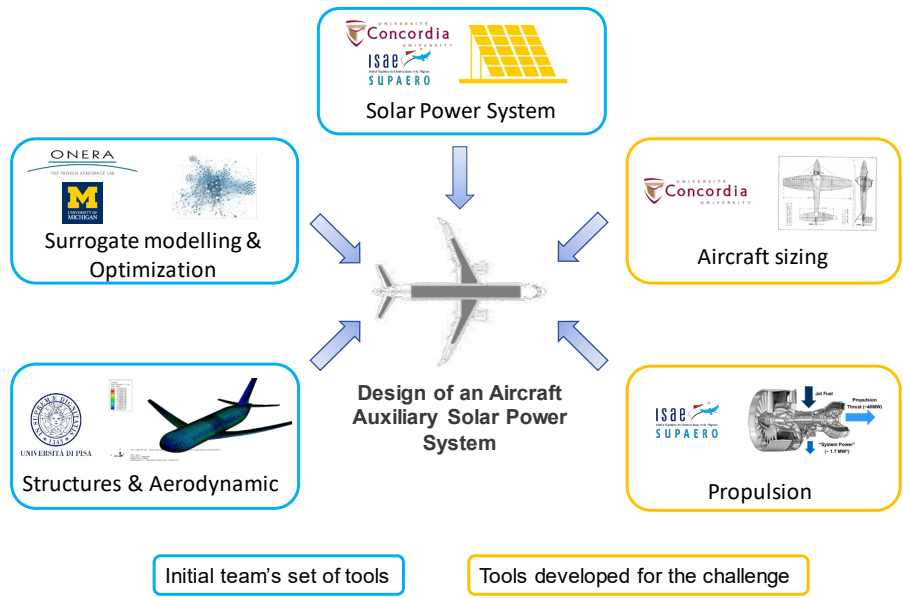

*Figure 4: Overview of tools used to solve the design problem*

Based on the tools, an MDO workflow has been established to evaluate the effect of implementing SPS on the baseline aircraft (Figure 5). The aircraft sizing tool represents the first element of the workflow. It provides the main characteristics of the aircraft to the other tools and evaluate the top-level aircraft parameters such as the Maximum Take-Off Weight (MTOW). The wing and fuselage characteristics are the inputs of the SPS tool. It determines the available solar power and the weight of the solar panels. These are the inputs of the propulsion tool and the structure tool respectively. The first one evaluates the amount of fuel that can be saved by using the solar power generated by the SPS to supplement the power offtakes during ground and cruise segments. In mean time, the weight of the SPS is passed to the aerodynamic and structure tools hosted by the University of Pisa. They use the flight conditions, the aircraft geometry, the SPS location and weight to compute an aero-structural analysis. In the meantime, an optimization of the empty weight is conducted. Finally, the fuel-weight savings and the empty weight are passed to the aircraft sizing tool that resizes the aircraft to maintain the same performance of the baseline aircraft.

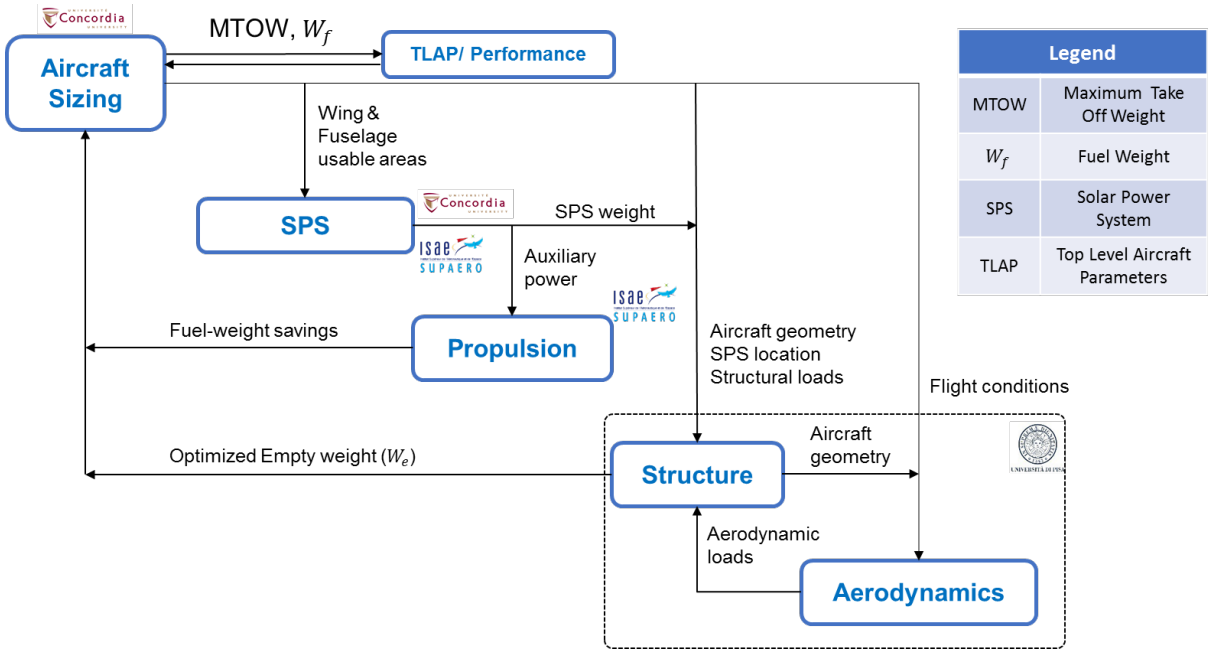

*Figure 5: MDAO workflow for the design problem*

The developed workflow defines the way the tools are connected and the inputs and outputs of every tools. The execution of the workflow requires common exchangeable inputs/outputs files to ensure the tool compatibility for each workflow iteration. It has been done using the CPACS files introduced in section 2.2.1. These files are exchanged between the tool's owner after tool execution with an updated version number for each iteration of the workflow.

### 3.3 Task B: Collaborative workflow development

Task B deals with workflow development support and team collaboration using AGILE advanced instruments for MDAO (KE-Chain, MDAO formalization, Visualization features). The KE-Chain platform was widely used from the outset of Task A. Documentation and requirements were specified on the platform and assigned to different team members to ensure compliance. The online platform aims to serve as an ultimate tool to configure the data model for the proposed project. Operating the KE-Chain platform as intended allows the operator to stay organized by uploading of user information which includes members involved in the project, assigning them roles within the project by allocating responsibility to each user to comply with requirements, upload competence data and a request log to identify missing files.

The first sections of the platform are very straightforward, mostly involves uploading user data, defining the different requirements and parameters of operating the tools, defining the tools used in the project; description, input, output information, user responsible for tool operation and data upload.

The second section is where the data model begins to be defined. As the AGILE project is being orchestrated using the common language of CPACS, KE-chain is designed to read these CPACS specific .xml files from each design tool and merge them as a single baseline file. This helps create a database of the CPACS files generated by each tool and help the user to visualize the combined hierarchy of data. The visualization is provided by the onboard VITMOS tool, which offers a nodal, hierarchy tree as shown in Figure 6.

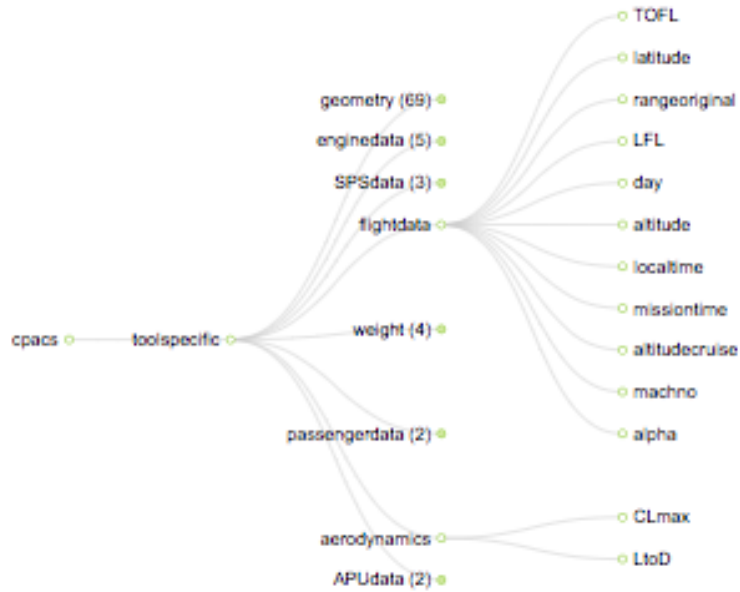

*Figure 6: Hierarchy tree view for the CPACS file of a tool*

A more useful function, provided next, is the function to match the earlier defined competences or tools with their appropriate input and output CPACS files. The user also has the option to import a CMDOWS file instead. By using the merge function, KE-chain's KADMOS script produces a CMDOWS file which defines the relationship between variables and tools. This is visually represented by a Repository Competence Graph (RCG), as shown in Figure 7.

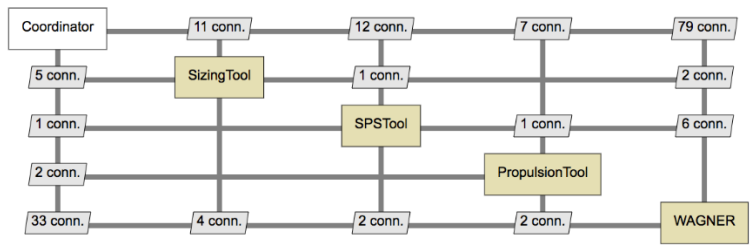

*Figure 7: RCG representation of the data model*

Combining the CPACS files and the tool definition, RCG helps map the different input and output variables coded in the CPACS files and how they interact between the different tools, where the top row blocks represent the input variables entering the tools and the left column blocks represent the output variable exiting the tools. A few more iterations are performed to produce a Fundamental Problem Graph (FPG), which streamlines and unclutters the data work flow by excluding, consolidating and addressing competence collisions (diagonal elements), and by selecting the key, desired design and state variable to be represented in the data work flow (off-diagonal elements in ). The FPG can be inspected to check the variable paths and competence information. Finally, the user defines the MDAO architecture settings to further optimize the data work flow, by selecting the basic solution strategy and decomposition of the coupling to be performed. The resulting work flow is shown in Figure 9. The final two sections involve the conversion of the CMDOWS files to be downloaded and used in the local RCE tool, and inspection of the downloaded files.

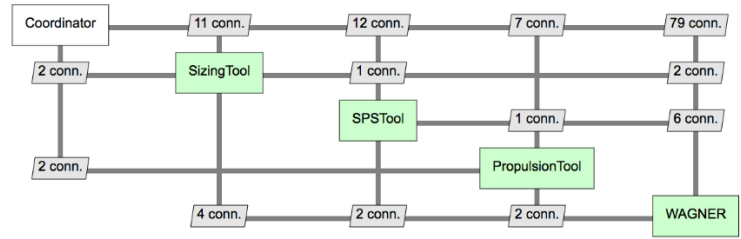

Figure 8: FPG representation of the data model

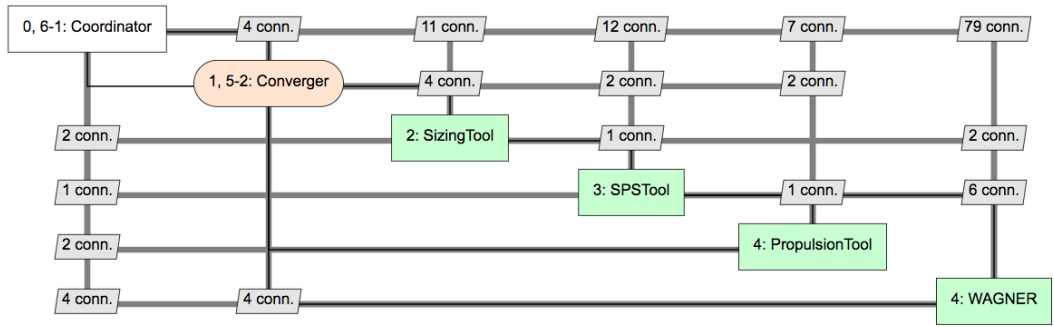

Figure 9: FPG with MDAO architecture applied - converged MDA

The final step of Task B dealt with the workflow tools integration into RCE using CPACS files as the exchange format. Tools were first tested by connecting them to each other to create simple workflows to test the platform. Then, a remote workflow was executed to complete the Task C.

### 3.4  Task C: MDAO using surrogate models

Task C of the Agile Challenge deals with the creation and use of surrogate models during the MDAO procedure. The aim of surrogate modeling is to create an analytical approximation of a model to reduce cost and computational time. To train surrogate models, the user must provide some inputs and outputs, called training points, based on Design of Experiments (DoE) techniques which are evaluated using the high-fidelity tools. The surrogate modelling and optimization procedure was supported by ONERA and University of Michigan team members through the use of the SEGOMOE optimizer [5].

The optimization objective was to maximize the fuel weight saved of the aircraft using solar panels to generate auxiliary power. The design variables of the optimization process are the solar panel efficiency, the wing area, the fuselage length, the fuselage diameter and the semi wingspan and tail span. The following table provides the upper and lower bounds of the design variables.

Table 2: Upper and lower bounds of the design variables

| Design variables | Lower Bounds | Upper Bounds | Units |
|---|---|---|---|
| Solar Panel Efficiency | 0.25 | 0.55 | - |
| Wing Area | 100 | 250 | m² |
| Fuselage Length | 32 | 40 | m |
| Fuselage Diameter | 3 | 4.5 | m |
| Semi wing span | 14 | 22 | m |
| Semi tail span | 5 | 8 | m |

This optimization was conducted with several constraints managed by SEGOMOE. The fuselage length to diameter and aspect ratio of the wing are controlled within defined ranges. The Maximum TakeOff Weight (MTOW) was evaluated and constrained to stay under the baseline value. Table 5 gives the ranges of these constraints. Other constraints are managed by each discipline within the specific tools such as maximum admissible Von Mises stress and wing tip displacement for structural analysis and the usable areas for solar panels installation within the SPS tool.

Table 3: Constraints specifications

| Constraints | Lower Bounds | Upper Bounds | Units |
|---|---|---|---|
| Fuselage length / Fuselage Diameter | 8 | 10.5 | - |
| Wing aspect ratio $\left(AR = \dfrac{Wing\ span^2}{Wing\ area}\right)$ | 6 | 15 | - |
| MTOW | - | 67585 | kg |

The optimization procedure was conducted in two steps. First, an initial DoE was generated to initialize the process and to conduct a sensitivity analysis. Then, based on an initial set of data, the optimizer defined the configurations to be evaluated. The workflow was then executed through the exchange and processing of CPACS files by each tool owner and the initial DoE generated by the SEGOMOE tool was followed. Figure 14 shows the results of the sensitivity analysis conducted on the initial DoE data. It helps to compare the influence of each design variable on the amount of fuel weight saved. The results show that the solar panel efficiency is the most influential variable since it changes the available solar power used to reduce the fuel consumption. This is consistent with the fact that solar power technology is the limiting factor for implementing a solar power system on an aircraft. Moreover, the results show also that the wing area and fuselage characteristics such as length and diameter have the same influence on the fuel weight saved. The conclusion of this preliminary analysis is that two strategies can be followed for the reduction of fuel weight using solar panels i.e. the optimization of the aircraft geometry or optimization of the solar panel efficiency.

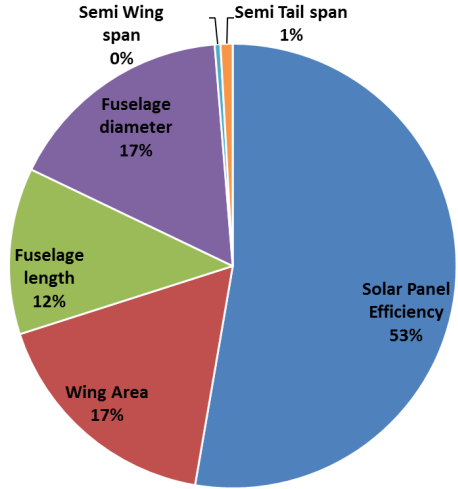

*Figure 10: Fuel weight saved sensitivity analysis conducted on the initial DoE data set*

Several iterations of the workflow were conducted to finalize the construction of surrogate models. The surrogate models were used to find an optimal aircraft according to optimization problem and constraints. Table 4 shows the baseline aircraft characteristics and the designed aircraft auxiliary solar power system. The resulting aircraft configuration has a 13% larger wing area and fuselage diameter than the baseline aircraft. The tail span is also larger by 21%. This has allowed 822 kg of fuel to be saved even if the solar power system results in a weight penalty. Although the results are promising, the tools used for the aircraft sizing were basic and more detailed analyses have to be conducted to reliably assess the potential of implementing solar panels on a commercial aircraft.

*Table 4: Optimization results and comparison with the baseline aircraft*

| Design variables | Baseline aircraft | Designed aircraft | Units |
|---|---|---|---|
| **Solar Panel Efficiency** | - | 0.41 | - |
| **Wing Area** | 113 | 128.7 | $m^2$ |
| **Fuselage Length** | 38 | 38.5 | m |
| **Fuselage Diameter** | 3.7 | 4.16 | m |
| **Semi wing span** | 17.5 | 17.1 | m |
| **Semi tail span** | 6.5 | 7.9 | m |
| **MTOW** | 67585 | 67585 | kg |
| **Fuel weight saved (SPS weight included)** | - | 822 | kg |

## 4 Discussions and lessons learned

The main objective of this AGILE challenge was to use the provided AGILE tools to design an aircraft auxiliary solar power system following three design tasks.

During the task A, an MDO workflow was developed based on the identification of design competencies of all collaborators involved in the project. Available tools were evaluated, and a design problem was formulated to suit the competences of the collaborating members. Tool gaps were identified in the workflow and a requirement for tool development was initiated.

Then, the tool development activity was tracked to completion on the KE-Chain platform during the task B of the challenge. The platform was also used to manage the implementation of the workflow by condensing all the individual tool into a single merged CPACS file.

Task C dealt with surrogate modeling and optimization procedures. The data generated by the workflow was incorporated into the creation of a surrogate model for the optimization activity. This optimization exercise is carried out using the SEGOMOE methodology and surrogate models developed for the workflow. The MDAO workflow execution data shows promising trends in terms of minimizing the fuel weight by implementing a solar power system on the AGILE Baseline Aircraft.

Several important lessons dealing with operating in a collaborative multidisciplinary environment were learnt throughout this challenge. First, the identification of competencies and the development of an MDAO workflow are time intensive. In addition, CPACS exchange files harmonization is recommended at early stages based on experience gained during the AGILE Challenge. Early definition of custom CPACS tag hierarchy and definition of tool inputs will make the overall process more efficient. The KE-Chain tool package helps better orchestrate the process, especially RCE has potential to make the optimization process faster.

Over the course of this project, the importance of team communication, competence identification and collaborative problem solving in design activities have been encountered and explored. This project has provided the authors with deeper insight and a greater appreciation of collaborative MDO tasks among heterogeneous teams of experts.

## 5   Acknowledgements

The authors are grateful to the AGILE consortium members and the other members of the team 3 who participated in the presented work: Marco Picchi Scardaoni, Karim Abu Salem, Mohamed Boulel, Pablo Rodrìguez.

## 6   References

[1]   P. D. Ciampa and B. Nagel, "AGILE the Next Generation of Collaborative MDO: Achievements and Open Challenges," in *AIAA Aviation Forum*, 2018.
[2]   E. Baalbergen, W. Lammen, N. Noskov, P.-D. Ciampa, and E. Moerland, "Integrated collaboration capabilities for competitive aircraft design," *MATEC Web Conf.*, vol. 233, p. 00015, Nov. 2018.
[3]   B. Aigner, I. van Gent, G. La Rocca, E. Stumpf, and L. L. M. Veldhuis, "Graph-based algorithms and data-driven documents for formulation and visualization of large MDO systems," *CEAS Aeronaut. J.*, vol. 9, no. 4, pp. 695–709, 2018.
[4]   P. Della Vecchia, B. Aigner, I. Van Gent, and P. D. Ciampa, "Collaborative Open Source Aircraft Design Framework for Education - AGILE Academy Initiative and Results," in *International Council of the Aeronautical Sciences*, 2018.
[5]   N. Bartoli *et al.*, "An adaptive optimization strategy based on mixture of experts for wing aerodynamic design optimization," in *18th AIAA/ISSMO Multidisciplinary Analysis and Optimization Conference*, 2017, pp. 1–17.